\documentclass[11pt]{amsart}
\usepackage[hmargin=2.5cm,vmargin=2.5cm]{geometry}
\usepackage{amsfonts}
\usepackage{amsthm}
\usepackage{amsmath}
\usepackage{amssymb}
\usepackage [T1]{fontenc}
\newtheorem{Theorem}{Theorem}[section]
\newtheorem{Definition}[Theorem]{Definition}
\newtheorem{Remark}[Theorem]{Remark}
\newtheorem{Example}[Theorem]{Example}
\newtheorem{lem}[Theorem]{Lemma}
\newtheorem{Proposition}[Theorem]{Proposition}
\newtheorem{Corollary}[Theorem]{Corollary}

\begin{document}
\title{Noncommutative Gr\"obner bases over  rings}
\author{\textbf{Andr\'e Mialebama Bouesso and Djiby Sow }}
\maketitle

\begin{center}\address{ D\'epartement de Math\'ematiques et Informatique \\ Laboratoire d'Alg\`ebre, de Cryptologie, de G\'eom\'etrie Alg\'ebrique et Applications\\ Universit\'e Cheikh Anta Diop\\
 BP 5005 Dakar Fann, S\'en\'egal}\\ E-mails :   \ miales2001@yahoo.fr, \  \ sowdjibab@yahoo.fr    \end{center}

\begin{abstract}

In this work,  it is proposed a method for computing Noncommutative Gr\"obner bases over a valuation n{\oe}therian ring. We have generalized the fundamental theorem on normal forms over an arbitrary ring. The classical method of dynamical commutative Gr\"obner bases is generalized for Buchberger's algorithm over  $R=\mathcal{V}\langle x_1,\ldots,x_m\rangle$ a free associative algebra with non-commuting variables, where $\mathcal{V}=\mathbb{Z}/n\mathbb{Z}$ or $\mathcal{V}=\mathbb{Z}$.

The process proposed,  generalizes previous known technics for the computation of Commutative Gr\"obner bases over a valuation n{\oe}therian ring and/or  Noncommutative Gr\"obner bases over a field.

\end{abstract}

\tableofcontents

\section*{\textbf{Introduction}}

Gr\"obner Bases is an algebraic technic that provides algorithmic solutions to a variety of problems in  Algebra and Algebraic Geometry. Noncommutative and Commutative Gr\"obner bases over a field provides many applications in computable algebras such as coding theory and cryptography.

 Many authors have  introduced and/or generalized Commutative Gr\"obner bases in different ways over fields or rings with zeros divisors (and non invertible elements more generally): see \cite{Buch}, \cite{Buch1}, \cite{Buch2A}, \cite{Buch2}, \cite{Cox}, \cite{HY},  \cite{Kapur},  and \cite{Y2}.

 Noncommutative G\"obner bases was also studied and developed by many authors (see: \cite{Cojocaru}, \cite{Green},  \cite{Green0},  \cite{Mora1} and \cite{Mora}).

 Most of concepts on noncommutative Gr\"obner bases over a field is analogous to the commutative case over a field. Nevertheless, some of the  important differences are:

 -  most ideals of noncommutative algebras do not have finite Gr\"obner bases,

 - it is possible to find a noncommutative Gr\"obner bases nonempty through a single polynomial.

 -  in some cases, the computation of the overlap relation (S-polynomial) of two polynomials  in Buchberger algorithm, is not possible.

 Throughout this paper,  we  propose a method for computing Noncommutative Gr\"obner bases over a valuation n{\oe}therian ring $\mathbb{Z}/n\mathbb{Z}$ and $\mathbb{Z}$.
 Our process generalizes previous known technics for the computation of Buchberger's algorithm in Commutative Gr\"obner bases over a valuation n{\oe}therian ring or in Commutative dynamical Gr\"obner bases over a principal ideal ring \cite{Y}, \cite{Y2} and in  Noncommutative Gr\"obner bases over a field \cite{Green0},\cite{Green}, \cite{Mora1}, \cite{Mora}. In this paper, we work with  $R=\mathcal{V}\langle x_1,\ldots,x_m\rangle$ a free associative algebra with non-commuting variables over a   ring $\mathcal{V}$.  Although most of the results that we present here hold for a wider class of noncommutative algebras.

This paper is structured as follows:

Section 1: Preliminaries: we adapt to  ring, the classical notions needed for Buchberger's algorithm in $R$ .

Section 2: We  introduce Noncommutative Gr\"obner bases and Reduced Gr\"obner bases  and we give also the fundamental theorem on normal forms over $R$.

Section 3: We give the division algorithm and Buchberger's algorithm for Noncommutative Gr\"obner bases over a valuation n{\oe}therian ring $\mathcal{V}$.

Section 4: We have generalized  Noncommutative Gr\"obner bases over  $\mathcal{V}=\dfrac{\mathbb{Z}}{n\mathbb{Z}}$.

Section 5: We adapt dynamical commutative Gr\"obner bases to Noncommutative Gr\"obner bases over the integers  $\mathcal{V}=n\mathbb{Z}$.

\section{\textbf{Preliminaries}}

In this section, we give some notions and notations that we will use in the sequel.

Let $\mathcal{V}$ be a commutative ring. $\mathcal{V}$ is said to be a valuation ring if for all $a,b\in \mathcal{V}$, $a$ divides $b$ or $b$ divides $a$.

A finite set of symbol is called alphabet. A finite sequence of elements of an alphabet $\sum $ is called a word. By a monomial, we mean a finite noncommutative word in the alphabet $\{x_1,\ldots,x_m\}$. We use the letter $\mathbb{M}$ to denote the set of monomials. We define the multiplication in the set $\mathbb{M}$ by concatenation. Let $p\in \mathbb{M}$, we denote by $Lth(p)$ the length of $p$ i.e the number of letter of the word $p$. Note that $Lth(w)=0$ if $w$ is an empty word.

Let $R=\mathcal{V}\langle x_1,\ldots,x_n\rangle$ be the free associative algebra with non-commuting variables defined over a ring $\mathcal{V}$, then $f\in R\Leftrightarrow f=\displaystyle\sum_{\alpha}a_\alpha p_\alpha$ as a finite sum,  where $p_\alpha\in \mathbb{M}$ and $a_\alpha \in \mathcal{V}$ with $p_\alpha \neq p_\alpha'$ if $\alpha \neq \alpha'$.
By a term $t$, we mean $t=ap$ where $a\in \mathcal{V}$ and $p\in \mathbb{M}$. We use the letter $\mathbb{T}$ to denote the set of terms. We define the multiplication in the set $\mathbb{T}$ as follows: if $t=ap$, $t'=a'p'$, with  $a, a'\in \mathcal{V}$ and $p, p'\in \mathbb{M}$, then $tt'=aa'pp'$ where $pp'$ is a concatenation of $p$ and $p'$ ($p$ and $p'$ do not commute). Let $t\in \mathbb{T}$. if $t=ap$ where $a\in \mathcal{V}$ and $p\in \mathbb{M}$, then $Lth(t)=Lth(p)$.

 A subset $I$ of $R$ is said to be a two-sided ideal of $R$ if:

 - $g_1-g_2 \in I $ for all $g_1,g_2\in I$

 - $f.g.h\in I$ for all $g\in I$ and all $f,h\in R$.

Let $D\subset R$, we denote by $\langle D\rangle=\{\sum f.g.h/ f,h\in R, g\in D\}$ the ideal generated by $D$ and we denote also by
$\langle\langle D\rangle\rangle=\{\sum_{i} \alpha_{i}.g_{i} / \alpha_{i} \in \mathcal{V}, g_{i}\in D\}$ the $\mathcal{V}$-submodule of $R$ generating by $D$.

         - An ideal $I\subset R$ is said to be a term-ideal if it is generated by elements in $\mathbb{T}$ of the form $ap$ with $a=1$ or $a\in \mathcal{ V}\setminus \mathcal{V}^{*}$ and $p\in \mathbb{M}$ (where $\mathcal{V}^{*}$ is the group of invertible elements of $\mathcal{V}$).

 Let $p $ and $q$ be two monomials of $\mathbb{M}$. We say that $p$ divides $q$ if there exists two monomials $m$ and $m'$ of $\mathbb{M}$ such that $q=m.p.m'$. Let $p\in \mathbb{M}$, we say that $p$ occurs in $f\in \mathcal{V}\langle x_1,\ldots,x_m\rangle$ if the coefficient of $p$ in $f$ is not zero. Let $t=ap $ and $t'=a'p'$ two terms of $\mathbb{T}$. We say that $t$ divides $t'$, if $a$ divides $a'$ in $\mathcal{V}$ and $p$ divides $p'$ in $\mathbb{M}$.

The set $\mathbb{M}$ of monomials is a mono\"id with the concatenation as mono\"id law.

A well-order $<$ on $\mathbb{M}$ is said to be admissible, if it satisfies the following conditions: for
all $p, q, r, s\in \mathbb{M}$, non empty :
\begin{itemize}
\item if $p<q$ then $pr<qr$
\item if $p<q$ then $sp<sq$ and
\item if $p=qr$ then $p>q$ and $p>r$.
\end{itemize}

In other words, $<$ is an admissible order if it is a well-order which is compatible with the mono\"id structure. Note that the admissible order generalizes the notion of monomial order for the commutative case.

Let $p$ and $q$ be two monomials in the finite alphabet $\{x_1,x_2,\ldots,x_n\}$. Considering that $x_1<x_2<\ldots<x_n$, the left graded lexicographic order (grlex) is defined as follows: $p<_{grlex}q$ for the left graded lexicographic order if:
\begin{itemize}
\item $Lth(p)<Lth(q)$ or
\item if $Lth(p)=Lth(q)$,  we find the biggest common left subword $m$ such that $p=m.w_1 \ \text{and}\  p_2=m.w_2$ where $w_1<w_2$ i.e the first symbol $x_i$ for $w_1$ is smaller than the first symbol $x_j$ for $w_2$ i.e $x_i<x_j$.
\end{itemize}
 The left graded lexicographic order is an admissible order.

 Throughout this paper we assume that $<$ is an admissible order.

\begin{itemize}
\item Let $q\in \mathbb{M}$, we say that $q$ is the leading monomial of $f$ and we note $q=LM(f)$ if $q$ occurs in $f\in \mathcal{V}\langle x_1,\ldots,x_m\rangle$ and $p<q$ for all monomials $p$ occurring in $f$.
\item  The coefficient of the leading monomial of $f$ is denoted $LC(f)$: it called the leading coefficient in $f$. The  term  $LT(f)=LC(f)LM(f)$ is the leading term of $f$.
\item  Let $E \subset R$, be a non empty set; and we define the following sets:
\begin{itemize}
 \item $LM(E):=\{LM(f)/ f\in E\setminus\{0\}\}$.
 \item $LT(E):=\{LT(g)/g\in E\setminus \{0\}\}$.
 \item $NonLM(E):=\mathbb{M}\setminus LM(E)$.
 \item $NonLT(E):=\mathbb{T}\setminus LT(E)$.
\end{itemize}
\end{itemize}


\section{\textbf{Noncommutative Gröbner bases over a commutative ring $\mathcal{V}$}}

In this section $R=\mathcal{V}\langle X_1,\ldots,X_n\rangle$ is a free associative algebra with noncommutative variables over a commutative ring $\mathcal{V}$ and $<$ an admissible order on the set of all monomials $\mathbb{M}$.

\subsection{\textbf{Noncommutative Gröbner basis}}

\begin{lem}
\label{lemme1}
 \begin{enumerate}
 \item Let $I=$ be a term-ideal of $R$. Then a term $b_\beta.q_\beta \in I$ if and only if there exists $\alpha_0p_{\alpha_0}$ in the set of generators of $I$ which divides $b_\beta.q_\beta $.
     \item Let $K =\langle\langle a_\alpha p_\alpha \ /\ p_\alpha \in \mathbb{M},  \alpha \in \Omega \subseteq \mathbb{N}^n \rangle\rangle$ be a $\mathcal{V}$-submodule of $R$.
       A term $b_\beta.q_\beta \in K$ if and only if there exists $\alpha_0\in \Omega$ such that $p_{\alpha_0}=q_\beta$ and $a_{\alpha_0}$ divides $b_\beta$. .
     \end{enumerate}
\end{lem}
\textbf{Proof} Obvious
   \hfill$\square$

\begin{Definition}
A subset $G\subset I$ (where $I$ is an ideal of $R$) is said to be a (two.sided) noncommutative Gr\"obner basis for $I$ with respect to $<$ if $$\langle LT(G)\rangle = \langle LT(I)\rangle$$.
\end{Definition}
How to construct a noncommutative Gr\"obner basis? The classical method is Buchberger's algorithm. We will see how to design Buchberger's algorithm over a valuation n{\"oe}therian ring ( section 3) and over a principal ideal ring (section 4).

In the following, we generalize the fundamental theorem on normal forms over $R$.
\begin{Theorem}  Let $I$ be a $\mathcal{V}$-submodule of $R$, then, as a $\mathcal{V}$-module, we have: $$R=I\oplus \langle \langle \mathcal{V}.LM(I)\setminus \mathcal{V}.LT(I)\rangle\rangle \oplus \langle \langle NonLM(I)\rangle\rangle,$$ where $\langle \langle X \rangle \rangle$ is the $\mathcal{V}$-submodule of $R$ generated by $X\subseteq R$.
\end{Theorem}

\underline{NB}: If $LC(f)$ is invertible for every $f\in I$, then  $\mathcal{V}.LM(I)\setminus \mathcal{V}.LT(I)=\emptyset$, thus, over a field  $R=I\oplus \langle \langle NonLM(I)\rangle \rangle$ and we retrieve the classical known result for normal forms.

\textbf{Proof}

Put $A=I$, $B=\langle \langle \mathcal{V}.LM(I)\setminus \mathcal{V}.LT(I)\rangle\rangle$ and $C= \langle \langle NonLM(I)\rangle\rangle$. First observe that $A\cap B=B\cap C=A\cap C=0$; know, we are going to prove that $A\cap (B\oplus C)=0$, $B\cap (A\oplus C)=0$ and $C\cap (A\oplus B)=0$.
\begin{enumerate}
\item Let us prove that $A\cap (B\oplus C)=0, B\cap (A\oplus C)=0$ and $C\cap (A\oplus B)=0$.

$(a)$  Suppose that  $f\in A\cap (B\oplus C)$ and $f\neq 0$.  We have $f\in A=I\Rightarrow LM(f)\in LM(I) \ \text{and} \ LT(f)\in LT(I)$ \ (*).
On the other hand, we have: $f\in B\oplus C \Rightarrow  LT(f)\in \mathcal{V}.LM(I)\setminus \mathcal{V}.LT(I) \ \text{or}\  LM(f)\in  NonLM(I) $.
Thus, using (*) \ we have \  $LT(f)\in LT(I)\cap \big[\mathcal{V}.LM(I)\setminus \mathcal{V}.LT(I)\big]=\emptyset $ \ \text{or}\ $LM(f)\in LM(I)\cap NonLM(I)=\emptyset $ which is impossible. Hence $A\cap (B\oplus C)=0$.

$(b)$ Suppose that  $f\in B\cap (A\oplus C)$ and $f\neq 0$.

We have  $f\in  B \Rightarrow LT(f), \ LM(f) \in \mathcal{V}.LM(I)\setminus \mathcal{V}.LT(I)$,  thus $LT(f)\notin LT(I)$ and $LM(f)\in LM(I)$ \  (**). On the other hand, we have:\ $f\in A\oplus C \Rightarrow LT(f) \in LT(I) \ \text{or} \ LM(f) \in NonLM(I)  $, which contradicts (**).  Hence $B\cap (A\oplus C)=0$.

$(c)$ Suppose that  $f\in C\cap (A\oplus B)$ and $f\neq 0$.

We have $f\in C= \langle \langle NonLM(I)\rangle \rangle \Rightarrow LM(f)\in NonLM(I) \Rightarrow LM(f)\notin LM(I)$ \ (***).  On the other hand, we have $f\in A\oplus B \Rightarrow LM(f) \in LM(I) \ \text{or} \  LM(f) \in \mathcal{V}.LM(I)\setminus \mathcal{V}.LT(I)$ (thus $LM(f)\in LM(I)$),  which is contradicts (***). Hence $B\cap (A\oplus C)=0$.

$(2)$ Now, let us show that $R=A\oplus B\oplus C=I\oplus \langle \langle\mathcal{V}.LM(I)\setminus LT(I)\rangle \rangle \oplus \langle \langle NonLM(I)\rangle \rangle$.

 Suppose that there exists $v\in R$ such that $v\notin A \oplus B\oplus C$ and let us prove that this fact is impossible. Define the set $H=\{LM(v)\ / \ v\notin A\oplus B\oplus C\}$, then $H \neq \emptyset$. Let $v_{0}$ be a minimal element with the property $LM(v_{0})\in H$.

 $(a)$ Suppose that $LT(v_{0})\in C=\langle \langle NonLM(I)\rangle \rangle$. Let $v_{1}=v_{0}-LT(v_{0})$ then we have $LM(v_{1})<LM(v_{0})\ (a*)$. By minimality from $(a*)$, we have $v_{1}\in A\oplus B\oplus C$ and $v_{0}=v_{1}+LT(v_0)\in A\oplus B\oplus C$. Which is impossible by definition of $v_0$.

$(b)$ Suppose that $LT(v_0)\notin C=\langle \langle NonLM(I)\rangle \rangle$.

We have $LT(v_0)\notin C=\langle \langle NonLM(I)\rangle \rangle \Rightarrow LM(v_0)\notin NonLM(I)$, thus $LM(v_0)\in LM(I)$, therefore   $LT(v_0)\in \mathcal{V}.LM(I)$.

     - Suppose that $LT(v_0)\in \mathcal{V}.LM(I)\setminus \mathcal{V}.LT(I)$. Let $v_{2}=v_0-LT(v_{0})\Rightarrow LM(v_{2})<LM(v_0) \ (b*)$. By minimality from $(b*)$,  we have $v_{2}\in A\oplus B\oplus C$ and $v_0=v_{2}+LT(v_{0})\in A\oplus B\oplus C$, which is impossible.

   - Suppose that $LT(v_0)\in \mathcal{V}.LT(I)$, we deduce that there exists $w\in I$ such that $LT(v_0)=LT(w)$. Let $v_{3}=v_0-w\Rightarrow LM(v_{3})<LM(v_0)$ \ (c*). By minimality from $(c*)$,  we have $v_{3}\in A\oplus B\oplus C$ and $v_0=v_{3}+w\in A+B+C$, which is impossible.

 We conclude that $R=A\oplus B \oplus C=I\oplus \langle \langle \mathcal{V}.LM(I)\setminus \mathcal{V}.LT(I)\rangle \rangle \oplus \langle \langle NonLM(I)\rangle \rangle$ as desired.
\end{enumerate}

 \hfill$\square$

\begin{Remark}
- Every element $f\in R$ has a unique decomposition $f=f_1+\mathcal{N}_{I}(f)$ where $f_1\in I \ \text{and}\ \mathcal{N}_{I}(f)=f_2+f_3\in \langle \langle \mathcal{V}.LM(I)\setminus \mathcal{V}.LT(I)\rangle \rangle \oplus \langle \langle NonLM(I)\rangle \rangle$. $\mathcal{N}_{I}(f)$ is called the normal form of $f$ relatively to $I$.

- Now, a question arise: how to compute the normal form of a polynomial? \\ As in the commutative case, we are going to see that division  by a set of generator of the ideal $I$ does not solve the problem in general, but division by a noncommutative Gr\"obner basis  of the ideal $I$  solves the problem over a suitable ring.
\end{Remark}

\subsection{\textbf{Division algorithm}}.

Division algorithm is already known in the non commutative case over a field \cite{green}. In the following we adapt this algorithm over a  ring.

Let $f\in R$, given an ordered set $F=\{f_1,\ldots,f_s\} \subset R$,  we propose a method to divide $f$ by $F$ i.e. we find nonnegative integers $t_1,\ldots,t_s$ and $u_{ij}, v_{ij},r \in R$ for $1\leq i\leq s$ and $1\leq j\leq t_i$ such that:
\begin{enumerate}
\item $f=\displaystyle\sum_{i=1}^{s}\displaystyle\sum_{j=1}^{t_i}u_{ij}f_iv_{ij} + r$.
\item $LM(f)\geq LM(u_{ij}f_iv_{ij})$ for all i and j.
\item $LT(f_i)$ does not divides any term occurring in $r$ for $1\leq i\leq s$. We will call $r$ a reminder of the division by $F$.
\end{enumerate}

\vspace{0,2cm}

\textbf{\underline{Algorithm1}}\\
\textbf{INPUT}: $F=\{f_1,\ldots,f_s\}$ (Ordered), $f$ and an admissible order $<$.\\
\textbf{OUTPUT}: $t_1,\ldots,t_s \ \in \mathbb{N}, u_{it_i},v_{it_i},r \in R$ such that $f=\displaystyle\sum_{i=1}^{s}\displaystyle\sum_{j=1}^{t_i}u_{ij}f_iv_{ij} + r$

\textbf{INITIALIZATION}: $t_1=\cdots=t_s=0$, $u_{it_i}=v_{i_i}=r=\emptyset$ and $h:=f$.

Divoccur:= False

\textbf{WHILE} $h\neq 0$ and Divoccur:=false Do,

\textbf{IF} $LT(f_{i})/LT(h)$ \big( with $LM(h)=u.LM(f_i).v$  and $LC(f_{i})/LC(h)$  for $1\leq i\leq s$ and $u,v \in \mathbb{M}$\big), then

$u_{it_i}:=[\frac{LC(h)}{LC(f_i)}]u$

$v_{it_i}:=v$

DIVOCCUR:=True

$h:=h - u_{it_i}.f_i.v_{it_i}$

\textbf{IF} DIVOCCUR $:=$ False, then

$r:=r+LT(h)$

$h:=h-LT(h)$
\begin{Example}
Let $(\mathbb{Z}/16\mathbb{Z})\langle x,y\rangle$. Let us divide $f=4(xy)^2-2xy$ by $f_1=3yxy+x^2$ and $f_2=2yx-6y$ with $x>_{grlex}y$. Then $f=-4x.f_1+4x^3-2xy$.

- If we start the division by $f_2$ we get $f=2x.f_2.y+12xy^2-2xy$.

- If we start the division by $f_1$ we get $f=-4xf_1+4x^3-2xy$
\end{Example}

\begin{Example}
Let $R=\mathbb{Z}\langle x,y,z\rangle$ and $I=\langle f_1=5xy-x, f_2=3x^2-xz\rangle $ be an ideal of $R$ generated by $F=\{f_1,f_2\}$. Let $f=30zx^2yx \in R$ and $> $ be the (left) graded-lexicographic order on $\mathbb{M}$ with $x>y>z$. The division of $f$ by $F$ yields

- $f=10z.f_2.yx+10zxzyx$ if we start by $f_{2}$,

- $f=6zx.f_1.x+2z.f_2.x+2(zx)^2$ if we start by $f_{2}$.
\end{Example}

\textbf{Notation}

If $F=\{g_{1}, \ldots, g_{n}\}$ is an ordered set in $R$ and $f \in R$, we denote by $r=\overline{f}^{F}$ a remainder of $f$ under the division by $F$ .

Note that from the above examples, we see that:

-The result of the division algorithm  depend on the order on $F$,

-The division algorithm doesn't allow to answer  the "ideal memberships problem" because if the remainder of the division $f$ by $F$ is $r$, we doesn't know if $r=\mathcal{N}_{I}(f)$ is the normal form of $f$.

In order to solve this two important problems, one must make the division by a noncommutative Gröbner basis as we will see in the following theorem.

\begin{Theorem} Suppose that $G$ is a noncommutative Gr\"obner basis of an ideal $I$ of $R$. Let $f\in R$ and assume that $F=\{g_1,\ldots,g_n\}=\{g\in G\ / \  LM(g)\leq LM(f)\}$. If $\overline{f}^F=r$ ($r\neq 0$) then $r$ is independent of the order of $g_1,\ldots,g_n$ in $F$. In fact \ \ $r=\mathcal{N}_{I}(f)$.
\end{Theorem}

\textbf{Proof}\ Consider that $\overline{f}^F=r \ (r\neq 0)$, then $LM(r)\leq LM(f)$, since $\langle LT(G)\rangle=\langle LT(I)\rangle$, we see that, for each $g\in G$, $LT(g)$ does not divide any term occurring in $r$. Hence $r\in  \langle \langle \mathcal{V}.LM(I)\setminus \mathcal{V}.LT(I)\rangle \rangle \oplus \langle \langle NonLM(I)\rangle \rangle$. On the other hand,  $f=\displaystyle\sum_{i=1}^{n}\displaystyle\sum_{j=1}^{t_i}u_{ij}g_iv_{ij} + r$ (From the division algorithm) with $LT(g_i), \ \forall 1\leq i\leq n$ does not divide any term occurring in $r$.  It is clear that $\displaystyle\sum_{i=1}^{n}\displaystyle\sum_{j=1}^{t_i}u_{ij}g_iv_{ij}\in I \ \text{and}\ r\notin I $, that implies $r\in  \langle \langle \mathcal{V}.LM(I)\setminus \mathcal{V}.LT(I)\rangle \rangle \oplus \langle \langle NonLM(I)\rangle \rangle$ since the decomposition is unique, we have $r=\mathcal{N}_{I}(f)$.

\hfill$\square$

\begin{Corollary}(Ideal membership problem)   Suppose that $G$ is a noncommutative Gr\"obner basis of an ideal $I$ of $R$. Let $f\in R$, then $f\in I$ if and only if $\overline{f}^G=0$.
\end{Corollary}
\textbf{Proof} Follows from the above proposition and the fundamental theorem on normal forms.

\hfill$\square$

\subsection{\textbf{Reduced noncommutative Gr\"obner basis}}

Let $f$ be an ideal, a set of term $N$ is a minimal generating set of term if $N=\{ap\in B_{J}/ \ \forall \ bq\in B_{J}, (bq/ap\Rightarrow p=q \ \text{and} \  \langle a\rangle=\langle b\rangle)\}$
\begin{Definition} Let $I$ be an ideal of $R$ and assume that  the term-ideal $ \langle LT(I) \rangle$ has a unique minimal terms generating set $T$ . We say that the set $G_{I,T}$ is the reduced noncommutative Gr\"obner basis for $I$ if $G_{I,T}=\big\{t-\mathcal{N}_{I}(t)/ t\in T \ \text{and} \ \mathcal{N}_{I}(t)\in \langle \langle \mathcal{V}.LM(I)\setminus \mathcal{V}.LT(I)\rangle \rangle \oplus \langle \langle NonLM(I)\rangle \rangle \big\}$.
\end{Definition}

\begin{Remark}
 The fact of having a unique minimal generating set of terms will guarantee the existence of the reduced noncommutative Gr\"obner basis.
\end{Remark}

\begin{Theorem} Let $I$ be an ideal of $R$ and assume that the term-ideal $\langle LT(I) \rangle$ has a unique minimal generating set of terms $T_{I}$.
 Let $G_{I,T}$ be the reduced noncommutative Gr\"obner basis for the ideal $I$ of $R$. Then the following hold.
\begin{enumerate}
\item $G_{I,T}$ is a noncommutative Gr\"obner basis for $I$.
\item If $g\in G_{I,T}$, then $LC(g)=1$ or $LC(g)\in \mathcal{V}\setminus \mathcal{V}^{*}$ is irreducible.
\item If $g\in G_{I,T}$ then $g-LT(g)\in \langle \langle \mathcal{V}.LM(I)\setminus \mathcal{V}.LT(I)\rangle \rangle \oplus \langle \langle NonLM(I)\rangle \rangle$.
\item $LT(G_{I,T})$ is the minimal terms generating set  of $ \langle LT(I) \rangle$.
\end{enumerate}
\end{Theorem}

\textbf{Proof}
\begin{enumerate}
\item $G_{I,T}$ is a noncommutative Gr\"obner basis for $I$?\\
We have to show that $\langle LT(G_{I,T})\rangle=\langle LT(I)\rangle$. It is obvious that $G_{I,T}\subseteq I$ and thus $LT(G_{I,T})\subseteq LT(I)$.
Let $f\in I\Rightarrow LT(f)\in LT(I)\subset \langle LT(I)\rangle=\langle T\rangle$, then there exists $t\in T$ such that $t/LT(f)\ (a)$. Put $g=t-\mathcal{N}_{I}(t)$ then $LT(g)\in LT(I)$.  Recall that $\mathcal{N}_{I}(t)=t_{2}+t_{3}$ with $t_{2}\in \langle \langle \mathcal{V}.LM(I)\setminus \mathcal{V}.LT(I)\rangle \rangle $ and $t_{3} \in  \langle \langle NonLM(I)\rangle \rangle$ and $g=t-t_{2}-t_{3}$. Since $LT(g)\in LT(I)$ then $LT(g)=LT(t_{2})$ or $LT(g)=LT(t_{3})$ are impossible, then $LT(g)=t$. Therefore, from $(a)$, we have $LT(g)/LT(f)\Rightarrow LT(f)\in \langle LT(G_{I,T})\rangle$, thus $\langle LT(I)\rangle=\langle LT(G_{I,T})\rangle$. Hence $G_{I,T}$ is a noncommutative Gr\"obner basis for $I$.
\item Obvious.

 \item If $g\in G_{I,T}$ then $g-LT(g)\in \langle \langle \mathcal{V}.LM(I)\setminus \mathcal{V}.LT(I)\rangle \rangle \oplus \langle \langle NonLM(I)\rangle \rangle$?\\
Let $g\in G $ then $ g=t-\mathcal{N}_{I}(t)\ \ (\alpha)$. From above $(1)$ we have seen that $\mathcal{N}_{I}(t)\neq LT(g)$ and $LT(g)=t \ \ (\beta)$. Thus from $(\alpha) - (\beta)$ we find $g-LT(g)=-\mathcal{N}_{I}(t)\in \langle \langle \mathcal{V}.LM(I)\setminus \mathcal{V}.LT(I)\rangle \rangle \oplus \langle \langle NonLM(I)\rangle \rangle$ as desired.
\item $LT(G)$ is the minimal generating set of terms of $\langle LT(I)\rangle$?\\
Since $G$ is a noncommutative Gr\"obner basis for $I$, we have $T\subset LT(G_{I,T})$ as in the above remark. Let $g\in G_{I,T}$ then $ g=t-\mathcal{N}_{I}(t)$, from $(2)$ we have seen that $LT(g)=t\in T\Rightarrow LT(G_{I,T})\subseteq T$. Thus $T=LT(G_{I,T})$ is the minimal generating set of term of $\langle LT(I)\rangle$.
\end{enumerate}
\hfill$\square$

\begin{Remark}
The above theorem guaranties the existence of reduced Gr\"obner basis but doesn't provide a procedure to compute it, since, until now, we are not able to compute the normal form of a polynomial. In the next section, we will have the necessary tools to compute the normal form of a polynomial.
\end{Remark}

\section{\textbf{Noncommutative Gr\"obner bases  over n{\oe}therian valuation ring }}

In this section we will give a way to construct a finite and an infinite noncommutative Gr\"obner basis by using the overlap relations which generalize S-polynomials for the commutative Gr\"obner bases \cite{Cox}. We will also recall the definition of Gr\"obner basis in the noncommutative case over a field (see: \cite{green}) and adapt it for a valuation ring.

The previous theorem prove that  noncommutative Gr\"obner bases allows to solve the "Ideal membership problem". But, how to compute a  Gr\"obner basis? For this, we need Buchberger's Algorithm.

In this section $R=\mathcal{V}\langle x_1,\ldots,x_n\rangle$ is a free associative algebra with non-commuting variables over a commutative ring $\mathcal{V}$ and $<$ an admissible order on the set of all monomials $\mathbb{M}$.

\textbf{Minimal generating sets of terms over n{\oe}therian ring }

\begin{Proposition}
Let $\mathcal{V}$ be a n{\oe}therian ring and $R=\mathcal{V}\langle x_1\ldots\ x_n\rangle$. Let $\leq$ be an admissible order on $\mathbb{M}$. If $J \subset R$  is a term-ideal of $R$,   then $J$ has a unique minimal generating set of terms which is $N=\{ap\in B_{J}/ \ \forall \ bq\in B_{J}, (bq/ap\Rightarrow p=q \ \text{and} \  \langle a\rangle=\langle b\rangle)\}$, where $B_{J}$ is the set of all terms in $J$ and  $\langle \lambda \rangle$ is the ideal of $\mathcal{V}$ generated by $\lambda \in \mathcal{V}$.

\end{Proposition}
\textbf{Proof} Let $S$ be the set of terms which generates $J$, i.e $J=\langle S\rangle$. Let $N=\{ap\in B_{J}/ \ \forall \ bq\in B_{J}, (bq/ap\Rightarrow p=q \ \text{and} \  \langle a\rangle=\langle b\rangle)\}$, where $\langle \lambda \rangle$ is the ideal of $\mathcal{V}$ generated by $\lambda \in \mathcal{V}$.

(1) \underline{Let us prove that $N$ is nonempty}.

 Put $\Sigma=\{p\in \mathbb{M} / \ \exists \ a\in  \mathcal{V}, \ ap\in S\}$, since $\leq $ is an admissible order on $\mathbb{M}$,   then $\Sigma$ has a minimal element $r$.

  - If $r\in J$ then $r\in B$, hence $r\in N$ thus $N$ is nonempty.

  - If $r\notin I$ then there exists $a\in \mathcal{V}\setminus \mathcal{V}^{*}$ (i.e $a$ is non-invertible) such that $ar\in S$ because $r\in \Sigma$. Thus, $C_{J}=\{\langle c\rangle, c \in \mathcal{V}\setminus \mathcal{V}^{*} / \ cr\in S\} \neq \emptyset$. Let $\langle c_{0}\rangle \subset \langle c_{1}\rangle  \subset \langle c_{2}\rangle \subset \ldots \subset  \langle c_{n}\rangle \subset \langle c_{n+1}\rangle \subset \ldots $ be an ascending chain of elements of $C_{J} $ since $\mathcal{V}$ is a n{\"oe}therian ring then  this chain is stationary i.e there exists $c_{i_{0}}$, such that $c_{i_{0}}=c_{i}$ for every $i\geq i_{0}$. Therefore  $c_{i_{0}}r$ is a minimal element and then  $c_{0}r$ belong to $N$ i.e $N$ is nonempty.

  (2) \underline{Let us prove that $N$ generates $J$}

Suppose that there exists $ap\in J$ such that $ap\notin \langle N\rangle$; by minimality there exists $bq\notin N$ such that $bq/ ap$ and $(p\neq q \ \text{or} \ \langle a\rangle \neq \langle b\rangle)$. But $bq / ap \Rightarrow (q\leq p \ \text{or} \ b/a)\Rightarrow (q\leq p \ \text{or} \  \langle a\rangle\subseteq \langle b\rangle )$ hence using the fact that $p\neq q$ or $\langle a\rangle \neq \langle b\rangle$ we have:

 $ap\notin \langle N \rangle \Rightarrow \ \exists \ bq\notin N$ such that $q<p$ or $\langle a\rangle \subsetneq \langle b\rangle$. \\
Now starting by $a_{1}p_{1}\notin \langle N \rangle$ and applying the above result recursively we have an infinite sequence of element $a_ip_i\notin \langle N\rangle$ such that $$\ldots /a_ip_i /\ldots /a_2p_2/a_1p_1$$ and at least one of the sequence $$\ldots\leq p_i\leq \ldots p_2\leq p_1 \ \quad (S_1)$$ and $$\langle a_1\rangle\subseteq \langle a_2\rangle\subseteq \ldots \subseteq \langle a_i\rangle\subseteq \ldots \ \quad (S_2)$$ is infinite.

 But the sequence $(S_1)$ is infinite, is impossible because $\leq $ is an admissible order,

  and also the sequence $(S_2)$ is infinite, is impossible because $\mathcal{V}$ is n{\oe}therian.

   We can conclude that $\forall \ ap\in B, \ ap\in \langle N\rangle$. Thus $N$ is a term generator set of $J$.

(3)\underline{ Finally let us prove that $N$ is minimal}
Suppose that there exist another generating set $N'$ of $J$ such that $N'\subset N$. If $a_0n_0\in N$, $a_0n_0\in J=\langle N'\rangle$ then there exists $a'n'\in N'$ such that $a'n'$ divides $a_0n_0$ and by definition of $N$ we have $\langle a_0\rangle=\langle a'\rangle$ and $n_0=n'$\ thus $a_0n_0=\alpha'a'n'\in N'$with $\alpha'\in \mathcal{V}$. Hence $N=N'$.
\hfill$\square$

\begin{Example}
Let consider $<_{grlex}$ order and Let $I\subset R$ be an ideal of $R$ such that $$I=\langle 2xy^{2}, 5x^{2}y, 3x^{2}, 4yx^{2}y, 3xy^{2}, 25xy^{2}, 25x^{2}y, 75x^{2}, 75y^{2}, 10x^{2}, 375y^{2}\rangle$$
\begin{enumerate}
\item over  $R=\dfrac{\mathbb{Z}}{5^{4}\mathbb{Z}}\langle x,y\rangle$, the unique minimal generating set of terms is  $N=\{xy^{2},  x^{2},  75y^{2}\}$.

\item over $R=\mathbb{Z}\langle x,y\rangle$, the unique minimal generating set of terms is  $$N=\{3x^2, 10x^2,75y^2,2xy^2,3xy^2,25xy^2,5x^2y,4yx^2y\}$$.

\end{enumerate}
\end{Example}

\begin{Remark}
\begin{itemize}
\item The unique minimal generating set of terms given in the above proposition is independent of any particular admissible order.
\item The unique minimal generating set of terms is not necessary finite, this differs from the commutative case with Dickson Lemma. For example the ideal $J=\langle xy^{i}x, i\in\mathbb{ N}, i>1\rangle$ has an infinite minimal generating set of terms.
\end{itemize}
\end{Remark}

\textbf{Buchberger's Algorithm over a valuation n{\oe}therian ring}.

First, let us generalizes the well known technic of overlap relation for non-commuting multivariate polynomials over a field.

\begin{Definition}  (\textbf{Overlap relation})
Let $f,g\in R=\mathcal{V} \langle X_1,\ldots,X_n\rangle $ where $\mathcal{V}$ is a valuation ring and $R$ a free associative algebra with $n$ non-commuting variables. Let $<$ be an admissible order on $\mathbb{M}$. Suppose that there are two monomials $p \ \text{and} \  q$:
\begin{enumerate}
\item $LM(f).p = q.LM(g)$
\item $LM(f)$ does not divide $q$ and $LM(g)$ does not divide $p$.\\ Then the overlap relation of $f$ and $g$ by $p \  \text{and} \ q$ is given by:
\begin{itemize}
\item $O(f,g,p,q) = \frac{LC(g)}{LC(f)}.f.p - q.g$ if $LC(f)/LC(g)$.
\item $O(f,g,p,q) = f.p - \frac{LC(f)}{LC(g)}.q.g$ if $LC(g)/LC(f)$.
\end{itemize}
\end{enumerate}
\end{Definition}

\begin{Definition}
 A subset $D\subset R$ is said to be LM-reduced if for all distinct elements $f, g\in D$, $LM(f)$ does not divide $LM(g)$ and vice versa.
 \end{Definition}

NB:  Over a field, every generators set can be LM-reduced in a new generators set, but that is not the case for a valuation ring (see \cite{Green}).

We present now the noncommutative version of Buchberger's algorithm over a valuation ring. We began by "Termination Theorem" which is a generalization of Bergman Diamond Lemma (\cite{Bergam}, \cite{green}).

\begin{Theorem} \ Let $\mathcal{V}$ be a n{\oe}therian valuation ring and $R=\mathcal{V}\langle x_1,\ldots,x_n\rangle$ Suppose that $G$ is a set of LM-reduced elements of $\mathcal{V}\langle x_1,\ldots,x_n\rangle$ such that every overlap relation $\overline{O(g,g',p,q)}^G=0$ with $g,g' \in G$, then $G$ is a noncommutative Gr\"obner basis for the ideal $I=\langle G\rangle$.
\end{Theorem}
\textbf{Proof} This proof is an adaptation to n{\oe}therian valuation ring of the proof in appendix of \cite{green}. Since there are some minor errors:

-  in the proof in  \cite{green} page 49,  for noncommutative Gr\"obner bases over field,

- and in the computation of S-polynomial of a single polynomial in \cite{Y} for commutative Gr\"obner bases over a valuation ring, (see corrigendum \cite{Y2}),

 we give  our proof for the shake of completeness . Let $f\in I$ and assume that $LT(f)$ is not divisible by $LT(g)$ for any $g\in G$. We need to show that this fact is impossible.

Assuming that $G$ is a generating set for $I$, we have $$f=\sum_{i,j}\alpha_{ij}p_{ij}g_{i}q_{ij}= \sum_{i,j}\alpha_{ij}p_{ij}LT(g_{i})q_{ij}+\alpha_{ij}p_{ij}\big[g_{i}-LT(g_{i})\big]q_{ij}  \quad (\Lambda)$$ with $g_{i} \in G$, $p_{ij}, \ q_{ij} \in \mathbb{M}$ and $\alpha_{ij} \in \mathcal{V}\setminus\{0\}$. Remark that  $LM[p_{ij}LT(g_{i})q_{ij}]= p_{ij}LM(g_{i})q_{ij}$ and  $LM[p_{ij}LT(g_{i})q_{ij}] > LM[p_{ij}\big(g_{i}-LT(g_{i})\big)q_{ij}]$ for each $(i,j)$. Let $K=\bigcup_{(i,j)/\alpha_{ij}\neq 0}\big\{p_{ij}LM(g_{i})q_{ij}\big\}$, $m=\#K$ and write $K=\big\{p_{i_{0}j_{0}}LM(g_{i_{0}})q_{i_{0}j_{0}}, p_{i_{1}j_{1}}LM(g_{i_{1}})q_{i_{1}j_{1}}\ldots, p_{i_{m}j_{m}}LM(g_{i_{m}})q_{i_{m}j_{m}}\big\}$ with \\ $p_{i_{l}j_{l}}LM(g_{i_{l}})q_{i_{l}j_{l}}>p_{i_{l+1}j_{l+1}}LM(g_{i_{l+1}})q_{i_{l+1}j_{l+1}}$, $0\leq l\leq m-1$.

Let $C_{l}=\sum_{(i,j)/p_{ij}LM(g_{i})q_{ij}=p_{i_{l}j_{l}}LM(g_{i_{l}})q_{i_{l}j_{l}}}LC(g_{i})\alpha_{ij}$ and

 $\Gamma_{l}=\big\{(i,j) /p_{ij}LM(g_{i})q_{ij}=p_{i_{l}j_{l}}LM(g_{i_{l}})q_{i_{l}j_{l}}\big\}$.

  Since $\mathcal{V}$ is a valuation ring, there exists one of the $LC(g_{i})$ occurring in $C_{l}$ which divides all the others and we denote it by $LC(g_{i^{(l)}})$ and by $p_{i^{(l)}j^{(l)}}LM(g_{i^{(l)}})q_{i^{(l)}j^{(l)}}$  the corresponding monomial \big[note that necessarily $p_{i^{(l)}j^{(l)}}LM(g_{i^{(l)}})q_{i^{(l)}j^{(l)}}=p_{i_{l}j_{l}}LM(g_{i_{l}})q_{i_{l}j_{l}}$ and $LC(g_{i^{(l)}})$ divides $C_{l}$, but it is possible to have $LM(g_{i^{(l)}})\neq LM(g_{i_{l}})$\big].

  If $C_{0}\neq 0$ then necessary

  $LT(f)=C_{0}p_{i_{0}j_{0}}LM(g_{i_{0}})q_{i_{0}j_{0}}=\frac{C_{0}}{LC(g_{i^{(0)}})}LC(g_{i^{(0)}})p_{i^{(0)}j^{(0)}}LM(g_{i^{(0)}})q_{i^{(0)}j^{(0)}}$, thus

   $LT(g_{i^{(0)}})=LC(g_{i^{(0)}})LM(g_{i^{(0)}})/ LT(f)$ which contradicts the hypothesis on $LT(f)$. Hence $C_{0}=0$ which implies that $\#\Gamma_{0}\geq 2$.  Put $z=\max \big\{l\ /\ 0\leq l\leq m-1, \ \#\Gamma_{l}\geq 2\big\}$, and considering all ways of rewriting $f$ as in $(\Lambda)$, we put $\Gamma$ to be the $\Gamma_{z}$ with the smallest size. Let $p^{*}$ the monomial corresponding to $\Gamma$ then $\#\Gamma$ is the number of occurrences of $p^{*}$ in $(\Lambda)$. We see that, if we consider all ways of rewriting $f$ as in $(\Lambda)$ then $p^{*}$ is the monomial which has the minimal number of occurrences. Therefore, there exists $p,q,p',q' \in \mathbb{M}$ and $g,g'\in G$ such that $$p^{*}=pLM(g)q=p'LM(g')q'$$

 Now, let us study all possible cases.

  \underline{Case:1} \ Suppose that $p<p'$

   \underline{Case:1.1} \ Suppose that $q\leq q'$ this forces that $LM(g')$ to divide $LM(g)$, which contradicts the fact that $G$ is LM-reduced.

    \underline{Case:1.2} \ Suppose that $q> q'$.

  From $p<p'$ and $q< q'$ we can write $p'=p\sigma$ and $q=\rho q' $ \ (with $\rho,\sigma \in \mathcal{M}$) \ and using that fact that $pLM(g)q=p'LM(g')q'$, we find $LM(g)\rho=\sigma LM(g')$.

 \underline{Case:1.2.1} \ Suppose that $LM(g)$ does not divide $\sigma$ and $LM(g')$ does not divide $\rho$. Then there is an overlap of $LM(g)$ and $LM(g')$ from $p^{*}$.
The corresponding overlap is :
 \begin{itemize}
\item $O(g,g',\rho,\sigma) = \frac{LC(g')}{LC(g)}.g.\rho - \sigma.g'$ if $LC(g)/LC(g')$,
\item and $O(g,g',\rho,\sigma) = g.\rho - \frac{LC(g)}{LC(g')}.\sigma.g'$ if $LC(g')/LC(g)$.
\end{itemize}

 Therefore, we have:

  \begin{itemize}
\item $p'g'q'=\frac{LC(g')}{LC(g)}pgq-p.O(g,g',\rho,\sigma)q'$ if $LC(g)/LC(g')$,
\item and $pgq=p.O(g,g',\rho,\sigma)q'+\frac{LC(g)}{LC(g')}.p'g'q'$ if $LC(g')/LC(g)$.
\end{itemize}

By assumption, $\overline{O(g,g',\rho,\sigma)}^{G}=0$, thus we have $O(g,g',\rho,\sigma)=\sum_{i}\widehat{\alpha}_{i}\widehat{p}_{i}\widehat{g}_{i}\hat{q}_{i}$ with $\widehat{\alpha}_{i}\in \mathcal{V}\setminus \{0\}$, $\widehat{p}_{i},\widehat{q}_{i} \in \mathbb{M}$ and $\widehat{g}_{i}\in G$, such that $LM(\widehat{p}_{i}\widehat{g}_{i}\widehat{q}_{i})<LM(g)\rho=\sigma LM(g')$.

Since  $LT(LC(g')pgq)=LC(g')LC(g)p^{*}=LT(LC(g)p'g'q')$, rewriting $pgq$ and $p'g'q'$ in this way, we can combine their leading term in order to lower the number of occurrences of $p^{*}$, which contradicts the minimality of occurrences of $p^{*}$.

 \underline{Case:1.2.2} \ \ Suppose that $LM(g)$  divides $\sigma$ or $LM(g')$  divides $\rho$, then from $LM(g)\rho=\sigma LM(g')$ we see that $LM(g)$  divides $\sigma$ and $LM(g')$  divides $\rho$ ( hence there is no overlap of $g$ and $g'$ in  $p^{*}$ and  $p'>pLM(g)$ and $q>LM(g')q'$).

 Now, write $g=\sum_{i}\alpha_{i}p_{i}+LT(g)$ and $g'=\sum_{j}\alpha'_{j}p'_{j}+LT(g')$. We have:

 \begin{itemize}
\item $\frac{LC(g')}{LC(g)}.pgq - p'g'q'=\frac{LC(g')}{LC(g)}\sum_{i}\alpha_{i}pp_{i}q - \sum_{j}\alpha'_{j}p' p'_{j}q'$ if $LC(g)/LC(g')$,
\item and $ pgq - \frac{LC(g)}{LC(g')}.p'g'q'=\sum_{i}\alpha_{i}pp_{i}q - \frac{LC(g)}{LC(g')}\sum_{j}\alpha'_{j}p'p'_{j}q'$ if $LC(g')/LC(g)$.
\end{itemize}

  This implies that:

  \begin{itemize}
\item $p'g'q'=\frac{LC(g')}{LC(g)}pgq-\frac{LC(g')}{LC(g)}\sum_{i}\alpha_{i}pp_{i}q + \sum_{j}\alpha'_{j}p' p'_{j}q'$ if $LC(g)/LC(g')$,
\item and $pgq=\sum_{i}\alpha_{i}pp_{i}q - \frac{LC(g)}{LC(g')}\sum_{j}\alpha'_{j}p'p'_{j}q'+\frac{LC(g)}{LC(g')}.p'g'q'$ if $LC(g')/LC(g)$.
\end{itemize}

   We deduce from the previous formulas that rewriting $pgq$ and $p'g'q'$ in this way, we can combine their leading term in order to lower the number of occurrences of $p^{*}$, which contradicts the minimality of occurrences of $p^{*}$.

 \underline{Case:2} Suppose that  $p=p'$ and we deduce from $p^{*}$ that $LM(g)q=LM(g')q'$ thus $LM(g)/LM(g')$ or  $LM(g')/LM(g)$; which contradicts the fact $G$ is LM-reduced

\underline{Case:3} Suppose that $p>p'$: this case is similarly to Case1, by symmetry.

\hfill$\square$

\vspace{0,2cm}

We are now ready to give the method  to construct  a noncommutative Gr\"obner basis.
 The procedure is the same than the commutative version, the difference is when computing the overlap relation of two polynomials. Note that in the noncommutative case, the noncommutative Gr\"obner bases is not necessary finite.\\
Let $\mathcal{V}$ be a n{\oe}therian valuation ring and $R=\mathcal{V}\langle x_1,\ldots,x_n\rangle$ a free associative algebra over $\mathcal{V}$.

Given $f_1,f_2,\ldots,f_k \in R=\mathcal{V}\langle x_1,\ldots,x_n\rangle$, let $I=\langle f_1,\ldots,f_k\rangle$ be a finitely, the algorithm produces a sequence of elements $g_1,g_2,\ldots,$ where $g_i=f_i$ for

$1\leq i\leq k$, and for $i>k$, $g_i\in I$ such that $LT(g_i)\notin \langle LT(g_1),\ldots,LT(g_{i-1})\rangle$.

\underline{\textbf{Algorithm 2}}

\begin{itemize}
\item Input: $\{f_1,f_2,\ldots,f_k\}$ a set of LM-reduced elements,
\item Output: $\{g_1,g_2,\ldots,g_k,\ldots\}$ a noncommutative Gr\"obner basis for $I=\langle f_1,\ldots,f_k\rangle$ for the admissible order\\
for $1\leq i\leq k$, do\\
$g_i:=f_i$\\$G:=\{g_1,\ldots,g_k\}$\\Count: k\\Do\\$\mathcal{H}:=G$\\For each pair of elements $h,h' \in \mathcal{H}$ and each overlap relation of $h,h'$,\\Do\\If $\overline{O(h,h',p,q)}^\mathcal{H}=r$ and $r\neq 0$, do\\Count:= count+1\\$g_{count}={r}$\\$G:=G\cup\{g_{count}\}$\\While $G\neq \mathcal{H}$\\$G$ is a noncommutative Gr\"obner basis for $I=\langle f_1,\ldots,f_k\rangle$ for the admissible order.
\end{itemize}
\textbf{Proof}: See the above theorem.

\begin{Proposition}
If the terms ideal of $I$ has a finite set of monomial generators, then the above algorithm terminates in a finite numbers of steps and yields a finite Gr\"obner basis.
\end{Proposition}
\textbf{Proof}: Similar to the proof of \cite{green} page 21. This case is similarly to the case of noncommutative Gr\"obner basis over a field and the fact that the ring is a n{\"oe}therian valuation ring guaranties that this algorithm terminates.

\begin{Example}
Let $f=3xyx-2xy$ be a polynomial in $(\mathbb{Z}/4\mathbb{Z})\langle x,y\rangle$ with noncommuting variables. Consider the left graded-lexicographic order with $x>y$ and let us construct a Gr\"obner basis for the ideal $I=\langle f\rangle$ of $(\mathbb{Z}/4\mathbb{Z})\langle x,y\rangle$. Set $G=\{f\}$, we can easily see that $LM(f)=xyx$ and $LM(f).yx=xy.LM(f)$

then $O_1=O(f,yx,xy)=2(xy)^2-2xy^2x \displaystyle\rightarrow _{f}-2xy^2x=f_1$, we replace the previous $G$ by $G=\{f,f_1\}$.

Since $LM(f).y^2x=xy.LM(f_1)$ then $O_2=O(f,f_1,y^2x,xy)=0$.

Since $LM(f_1).yx=xy^2.LM(f)$ then $O_3=O(f_1,f,yx,xy^2)=0$.

Since $LM(f_1).y^2x=xy^2.LM(f_1)$ then $O_4=O(f_1,y^2x,xy^2)=0$.

Hence $G=\{3xyx-2xy, -2xy^2x\}$ is a Gr\"obner basis for $I$ in $(\mathbb{Z}/4\mathbb{Z})\langle x,y\rangle$.
\end{Example}

\begin{Example}
Let $R=\mathbb{Z}/9\mathbb{Z}\langle w,x,y,z\rangle$. Consider the left-graded-lexicographic order with $x>y>z>w$. Let us construct a Gr\"obner basis for the ideal $I=\langle f_1=3yzwx-2yx, f_2=4xy-zw\rangle$ of $R$. Set $g_1:=f_1, g_2:=f_2\ \text{and}\ G=\{g_1,g_2\}$. We can easily see that $LM(g_1)=yzwx,\ LM(g_2)=xy\ \text{and}\ LM(g_1).y=yzw.LM(g_2)$ then $O_1=O(g_1,g_2,y,yzw)=3y(zw)^2-2yxy \Rightarrow \overline{O(g_1,g_2,y,yzw)}^G=3y(zw)^2-4yzw=g_3$.

The previous $G$ is replaced by $G=\{g_1,g_2,g_3\}$.

 Since $LM(g_2).(zw)^2=x.LM(g_1)$ then $O(g_2,g_1,(zw)^2,x)=-3(zw)^3-4xyzw\Rightarrow \overline{O(g_2,g_1,zwx,x)}^G=-3(zw)^2-4zwx=g_4$. The previous $G$ is replaced by $G=\{g_1,g_2,g_3,g_4\}$. Since $\overline{g_3}^{g_4}=0$, then $g_3$ is removed from $G$ and the new $G$ is $G=\{g_1,g_2,g_4\}$.

 Since there is no more overlap relation in $G$ then $G$ is a Gr\"obner basis for $I$ in $R$.

\end{Example}

\section{\textbf{Noncommutative Gr\"obner bases over} $R=\frac{\mathbb{Z}}{n\mathbb{Z}}\langle x_1,\ldots,x_m\rangle$}

Let $\mathcal{V}_{i}$ for $1\leq i\leq n$ be a finite family of rings such that we are able to compute Noncommutative Gr\"obner basis $G_{i}$ in $R_{i}=\mathcal{V}_{i}\langle x_1,\ldots,x_m\rangle$ which solves the "Ideal membership problem" for $\langle G_{i}\rangle$, i.e $f_{i}\in \langle G_{i}\rangle \Leftrightarrow \overline{f}^{G_{i}}_{i}=0 $ for each $1\leq i\leq k$.

Define the canonical projection $\pi_{i}:\mathcal{V}=\displaystyle\prod^{k}_{j=1}\mathcal{V}_{j} \rightarrow \mathcal{V}_{i}: a=(a_{1},\ldots,a_{k})\mapsto \pi_{i}(a)=a_{i} $. Then this projection extends naturally in a projection $\pi_{i}:\prod^{k}_{j=1}R_{j} \rightarrow R_{i}$. If $J$ is an ideal of $\displaystyle\prod^{k}_{j=1}R_i$ then $J_{i}=\pi_{i}(J)$ is an ideal of $R_{i}$, and $f\in J \Leftrightarrow f_{i}=\pi_{i}(f) \in J_{i}$ for each $1\leq i\leq n$. We know that if $J_{i}=\langle G_{i}\rangle$, then $f\in J \Leftrightarrow \overline{f}^{G_{i}}_{i}=0, \ \forall 1\leq i\leq k$.

Now, since $G=\big(G_{1},\ldots, G_{k}\big)$ guaranties the "Ideal membership problem" for $J$, we call $G$ the Dynamical Noncommutative Gr\"obner basis for $J$.

If we have a ring isomorphism $\psi: \mathcal{A} \rightarrow \mathcal{V}=\prod^{k}_{j=1}\mathcal{V}_{j}$ then it extends naturally in isomorphism $\psi: \mathcal{A}\langle x_1,\ldots,x_m\rangle \rightarrow \mathcal{V}\langle x_1,\ldots,x_m\rangle$. Now let $I$ be an ideal of $\mathcal{A}\langle x_1,\ldots,x_m\rangle$, then $f\in I\Leftrightarrow \psi(f) \in J=\psi(I)$. In this case, we say also that $G=\big(G_{1},\ldots, G_{k}, \psi \big)$ is the $\psi $-Dynamical Noncommutative Gr\"obner basis for $I$ where $G_{i}$ is a Noncommutative Gr\"obner basis for $J_{i}=\pi_{i}\circ \psi(I)$ for each $1\leq i\leq k$.

\vspace{0,2cm}

\textbf{Applications in $\frac{\mathbb{Z}}{n\mathbb{Z}}$}

\vspace{0,2cm}

Let $n=p^{\alpha_{1}}_{1}\ldots p^{\alpha_{k}}_{k}$ be an integer (where $p_{k}$ is a prime and $\alpha_{j}\in \mathbb{N}$) and $\psi: \frac{\mathbb{Z}}{n\mathbb{Z}}\rightarrow \prod^{k}_{j=1}\mathbb{Z}/p^{\alpha_{j}}_{j}\mathbb{Z}$ be the classical isomorphism from Chinese remainder theorem. Since $\mathbb{Z}/p^{\alpha_{j}}_{j}\mathbb{Z}$ is a valuation n{\"oe}therian ring, then from section 3, we are able to compute Noncommutative Gr\"obner bases in   $\mathbb{Z}/p^{\alpha_{j}}_{j}\mathbb{Z}$. Therefore, we can solve the "Ideal membership problem"  for an ideal in $\frac{\mathbb{Z}}{n\mathbb{Z}}\langle X_1,\ldots,X_m\rangle $ by using dynamical process.

\vspace{0,2cm}

\textbf{Example}

\vspace{0,2cm}

Let $n=24=3.2^{3}$ and $I=\langle f_{1}=14xy^{2}x-16yx^{2}, f_{2}=22x^{2}y^{2}-36yx\rangle$ be the ideal of $\frac{\mathbb{Z}}{24\mathbb{Z}}\langle x,y\rangle$, fix the left graded-lexicographic order with $x>y$. By the Chinese remainder theorem we can define the following isomorphism: $\varphi: \mathbb{Z}/24\mathbb{Z} \rightarrow
(\mathbb{Z}/3\mathbb{Z})\times(\mathbb{Z}/8\mathbb{Z}): x \mod 24\mapsto (x \mod 3, x \mod 8),$

Moreover, we have: ${\varphi}^{-1}: (\mathbb{Z}/3\mathbb{Z})\times(\mathbb{Z}/8\mathbb{Z})\rightarrow \mathbb{Z}/24\mathbb{Z}: (x \mod 3, x \mod 8)\mapsto (16x+9y) mod 24.$

Let us denote by $I_1= {\pi}_1(I)=\langle h_1=2xy^2x-yx^2, h_2=x^2y^2\rangle \subset (\mathbb{Z}/3\mathbb{Z})\langle x,y \rangle$ and $I_2= {\pi}_2(I)=\langle h'_1=6xy^2x, h'_2=6x^2y^2+4yx\rangle \subset (\mathbb{Z}/8\mathbb{Z})\langle x,y \rangle$.

We can easily see that $G_1=\{2xy^2x-yx^2, x^2y^2, xyx^2\}$ is a Gr\"obner basis for $I_1$ in $(\mathbb{Z}/3\mathbb{Z})\langle x,y \rangle$ and $G_2=\{6xy^2x, 6x^2y^2+4yx, 4xy^3x, 4yx^2\}$ is a Gr\"obner basis for $I_2$ in $(\mathbb{Z}/8\mathbb{Z})\langle x,y \rangle$. Thus $G=\{G_1, G_2, \varphi\}$ is a $\psi $-dynamical Gr\"obner basis for $I$ in $(\mathbb{Z}/24\mathbb{Z})\langle x,y \rangle$.


\section{\textbf{Noncommutative Gr\"obner bases over the integers}}

\subsection{\textbf{Special and Dynamical noncommutative Gr\"obner bases}}

\begin{Definition}

 \begin{itemize}
        \item $S$ is a multiplicative subset of a ring  $\mathcal{V}$ if $S \subseteq \mathcal{V}$, $1\in S$ and $ \forall x, y\in S$, \  $xy \in S$.
\item  $\mathcal{M}(x_{1},\ldots, x_{r}) = \{x_{1}^{n_{1}}\times \ldots \times x_{r}^{n_{r}}, \ n_{i} \in \mathbb{N}\}$ is the multiplicative subset of $\mathcal{V}$ generated by $\{x_{1},\ldots , x_{r}\}$ where $x_{1},\ldots, x_{r} \in \mathcal{V}$. Briefly, we denote $x^{\mathbb{N}}=\mathcal{M}(x)=\{x^{n}, \ n\in \mathbb{N}\}$.
\item Let $S$ be a multiplicative subset of a ring  $\mathcal{V}$. The ring $S^{-1}\mathcal{V} = \{\frac{x}{s}, x\in \mathcal{V}; s\in S\}$ is the localization of $\mathcal{V}$ relatively to $S$.
\item Let $x \in \mathcal{V}$, we denote by  $\mathcal{V}_{[x]}$, the localization of $\mathcal{V}$ relatively to the multiplicative subset $x^{\mathbb{N}}$. Moreover, one can define by induction $\mathcal{V}_{[x_{1},x_{2},\ldots, x_{k}]} := (\mathcal{V}_{[x_{1},x_{2},\ldots, x_{k-1}]})_{[x_{k}]}$ for $x_{1},\ldots, x_{k} \in \mathcal{V}$.
\item The multiplicative subsets $S_{1},\ldots, S_{k}$ of a ring  $\mathcal{V}$ are  called comaximal if $ \forall \ s_{1}\in S_{1},\ldots, s_{n} \in S_{n}; \ \exists \ v_{1},\ldots, v_{n} \in \mathcal{V}$ such that $\sum_{i=1}^{n}v_{i}s_{i} = 1$.
 \end{itemize}
\end{Definition}

\begin{Definition}

  -  Let $\mathcal{V}$ be a Principal ideal ring, $f, g \in \mathcal{V}\langle x_1,\ldots, x_n\rangle \backslash\lbrace 0\rbrace$, $I=\langle f_1\ldots ,f_s\rangle$ be a nonzero finitely generated ideal of $\mathcal{V}\langle X_1,\ldots, X_n\rangle$, and $>$ an admissible order.

- If $G=\lbrace g_1,\ldots,  g_t\rbrace$  and $\lbrace LC(g_1),\ldots,  LC(g_t)\rbrace$ are totally ordered under division then it is possible to define for each $(i,j)$, the overlap of $g_i$ and $g_j$ such as on valuation rings. Therefore, we are able to generalize noncommutative Gr\"obner basis to Principal ideal rings.

    \begin{itemize}
        \item For $g_1,\ldots ,g_t\in \mathcal{V}\langle X_1,\ldots, X_n\rangle$, $G=\lbrace g_1,\ldots,  g_t\rbrace$ is said to be a special noncommutative Gr\"obner basis for $I$ if $I=\langle g_1,\ldots,  g_t\rangle$, the set $\lbrace LC(g_1),\ldots,  LC(g_t)\rbrace$ is totally ordered under division and for each $(g,g')\in G\times G$, $\overline{O(g,g',p,q)}^{G}=0$.
        \item A set $G=\lbrace (S_1,G_1),\ldots, (S_k,G_k)\rbrace$ is said to be a dynamical noncommutative Gr\"obner
         basis for $I$ if $S_1,\ldots, S_k$ are finite comaximal multiplicative subsets of $\mathcal{V}$ and in each localization $(S_{i}^{-1}\mathcal{V})[X_1,\ldots, X_n]$, $G_i$ is a special noncommutative Gr\"obner basis for $S^{-1}I$.
    \end{itemize}
  \textbf{NB} If $\mathcal{V}$ be a Principal ideal ring then each localization $S_{i}^{-1}\mathcal{V}$ is a Principal ideal ring
\end{Definition}

\begin{Proposition}
   Let $\mathcal{V}$ be a Principal ideal ring, $I=\langle f_1\ldots ,f_s\rangle$  a nonzero finitely-generated
    ideal of $\mathcal{V}\langle X_1,\ldots, X_n\rangle$. Let $f\in \mathcal{V}\langle X_1,\ldots, X_n\rangle$ and fix an admissible order. Suppose that $G=\lbrace g_1,\ldots,  g_t\rbrace$ is a special noncommutative Gr\"obner basis for $I$ in $\mathcal{V}\langle X_1,\ldots, X_n\rangle$. Then, $f\in I$ if and only if $\overline{f}^{G}=0$.

\end{Proposition}
\textbf{Proof} Similar to the one of Special commutative Gr\"obner basis over a Principal ideal ring in \cite{Y}  and in corrigendum \cite{Y2}.
\hfill$\square$

\begin{Theorem}
Let $\mathcal{V}$ be a Principal ideal rings, $I=\langle f_1\ldots ,f_s\rangle$  a nonzero finitely-generated ideal of $\mathcal{V}\langle x_1,\ldots, x_n\rangle$, $f\in \mathcal{V}\langle x_1,\ldots, x_n\rangle$ and fix an admissible order on $\mathcal{V}\langle x_1,\ldots, _n\rangle$. Suppose that $G=\lbrace (S_1,G_1),\ldots, (S_k,G_k)\rbrace$ is a dynamical noncommutative Gr\"obner basis for $I$ in $\mathcal{V}\langle x_1,\ldots, x_n\rangle$. Then, $f\in I$ if and only if $\overline{f}^{G_i}=0$ in $(S_{i}^{-1}\mathcal{V})[X_1,\ldots, X_n]$ for each $1\leq i\leq k$.

\end{Theorem}

\textbf{Proof}: Similar to the one of Dynamical commutative Gr\"obner basis over a Principal ideal ring in \cite{Y}  and in corrigendum \cite{Y2}.
\vspace{0.2cm}

\subsection{\textbf{ Buchberger's algorithm for Dynamical noncommutative Gr\"obner basis}}.

We are now ready to present the algorithm to construct a Special and a Dynamical noncommutative Gr\"obner basis for an ideal $I=\langle f_1,\ldots,f_s\rangle$ of $R=\mathbb{Z}\langle X_1,\ldots,X_n\rangle$.

Comparatively to n{\oe}therian valuation ring, we can find two incomparable (under division) elements $a, \ b \in \mathcal{V}$. In this case, one should compute $d = a \wedge b$, factorize $a = da'$, $b = db'$, with $a'\wedge b' = 1$, and then
open two branches from $\mathcal{V}$, the computations are pursued in $\mathcal{V}_{[a']}$ and $\mathcal{V}_{[b']}$.

\vspace{0,2cm}

\textbf{Case of overlap relation} The overlap relation of $f$ and $g$ by $p \  \text{and} \ q$ is given as follow:

- If $LC(f)$ and $LC(g)$ are  comparable in $\mathcal{V}$ under the division order then apply the classical definition (as in a valuation ring);

- If $LC(f)$ and $LC(g)$ are not comparable in $\mathcal{V}$ under the division order then:

\begin{enumerate}
\item write $LC(f)=(LC(f)\wedge LC(g)).a'$ and  $LC(g)=(LC(f)\wedge LC(g)).b'$ where $a'\wedge b'=1$
and  $LC(f)\wedge LC(g)=\gcd(LC(f), LC(g))$ and open two branch:

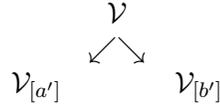
\begin{center}
$\mathcal{V} $ \\$\swarrow \searrow $ \\ $\mathcal{V}_{[a']} \hspace{1,5cm} \mathcal{V}_{[b']} $
\end{center}
\item in $\mathcal{V}_{[a']}=\{\frac{c}{a'^n}/ c\in \mathcal{V} \ \text{and}\ n\in \mathbb{N}\},\ a'$ is invertible and $LC(f)$ divides $LC(g)$, then the overlap relation is: \\
$O(f,g,p,q) = \frac{LC(g)}{LC(f)}.f.p - q.g$.
\item in $\mathcal{V}_{[b']}=\{\frac{c}{b'^n}/ c\in \mathcal{V} \ \text{and}\ n\in \mathbb{N}\},\ b'$ is invertible and $LC(g)$ divides $LC(f)$, then the overlap relation is: \\
$O(f,g,p,q) = f.p - \frac{LC(f)}{LC(g)}.q.g$
\end{enumerate}

\vspace{0,2cm}

\textbf{Case of division algorithm} For the division  of $f$ by $g$, if one has to divide $LT(f)=LC(f)LM(f)$ by $LT(g)=LC(g)LM(g)$ with $LM(g)/LM(f)$  and  $LC(f)$ and $LC(g)$ are not comparable in $\mathcal{V}$ under the division order then:

\begin{enumerate}
\item write $LM(f)=uLM(g)v$ with $u,v\in \mathbb{M}$,
\item write $LC(f)=(LC(f)\wedge LC(g)).a'$ and  $LC(g)=(LC(f)\wedge LC(g)).b'$ where $a'\wedge b'=1$
and  $LC(f)\wedge LC(g)=\gcd(LC(f), LC(g))$ and open two branch:

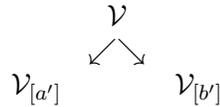
\begin{center}
$\mathcal{V} $ \\$\swarrow \searrow $ \\ $\mathcal{V}_{[a']} \hspace{1,5cm} \mathcal{V}_{[b']} $
\end{center}
\item in $\mathcal{V}_{[a']}$:\  $f = \frac{a'}{b'}u.g.v - r$ and the division is pursue were $f$ will be replaced by $r$.
\item in $\mathcal{V}_{[b']}$: \ $LT(f)$ is not divisible by $LT(g)$ and therefore $f=\overline{f}^{\{g\}}$
\end{enumerate}

  \vspace{0,2cm}
\begin{Example} \textbf{Overlap relation}

Let $f_1=6yzwx-2yx$ and $f_2=4xy-5zw$ be two polynomials in $R=\mathbb{Z}\langle x,y,z,w\rangle$ such that $w<x<y<z$, fix the left graded-lexicographic
order. We have: $LM(f_1).y=yzw.LM(f_2)$ and $LC(f_1)$ does not divides $LC(f_2)$ and vice versa, in this case we have:
$LC(f_1)=(LC(f_1)\wedge LC(f_2)).a'=(6\wedge 4).3$ and $LC(f_2)=(LC(f_1)\wedge LC(f_2)).b'=(6\wedge 4).2$ .
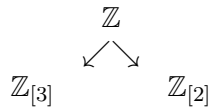
\begin{center}
$\mathbb{Z} $ \\$\swarrow \searrow $ \\ $\mathbb{Z} _{[3]} \hspace{1,5cm} \mathbb{Z} _{[2]} $
\end{center}

First case: In $\mathbb{Z}_{[3]}=\{\frac{a}{3^n}/ a\in  \mathbb{Z}, n\in \mathbb{N}\}$, $3$ is invertible and $LC(f_1)$ divides $LC(f_2)$ then
$$O(f_1,f_2,y,yzw)=\frac{LC(f_2)}{LC(f_1)}.f_1.y-(yzw).f_2=5yzwxy-\frac{4}{3}yxy$$

Second case:  In $\mathbb{Z}_{[2]}=\{\frac{a}{2^n}/ a\in \mathbb{Z}, n\in \mathbb{N}\}$, $2$ is invertible and $LC(f_2)$ divides $LC(f_1)$ then
$$O(f_1,f_2,y,yzw)=f_1.y-\frac{LC(f_1)}{LC(f_2)}.(yzw).f_2=\frac{15}{2}y(zw)^2-2yxy$$
\end{Example}

 \vspace{0,2cm}

\begin{Example}   \textbf{Buchberger's algorithm}
Let $R=\mathbb{Z}\langle x,y\rangle$ be a free associative algebra, $I=\langle f_1=6xyx-8xy, f_2=4xy-3yx\rangle$ be an ideal of $R$ and we fix an admissible order with $x>_{grlex}y$. Set $G=\{g_1=6xyx-8xy, g_2=4xy-3yx\}$, our goal is to construct a noncommutative Gr\"obner basis for $I$ in $R$. We have: $LM(g_1).yx=xy.LM(g_1)$, then $\overline{O(g_1,yx,xy)}^G=6yx^2y-y^2x=g_3$ and $G:=G\cup \{g_3\}=\{g_1,g_2,g_3\}$. Notice that $LC(g_1)=6,\ LC(g_2)=4$ and both are incomparable under division in $\mathbb{Z}$, then we open from $\mathbb{Z}$ two branches $\mathbb{Z}_{[2]}\ \text{and}\ \mathbb{Z}_{[3]}$ and pursued the computation of overlap relations in each branch:

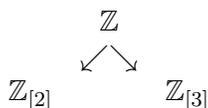
\begin{center}
$\mathbb{Z} $ \\$\swarrow \searrow $ \\ $\mathbb{Z} _{[2]} \hspace{1,5cm} \mathbb{Z} _{[3]} $
\end{center}

\begin{itemize}
\item \underline{In the ring $\mathbb{Z} _{[2]}$}, the set $G_1=\{6xyx-8xy, 4xy-3yx, \frac{-9}{2}yx^2+6yx, \frac{3}{2}y^2x\}$ is a Gr\"obner basis for $(2^{\mathbb{N}})^{-1}I$ in $((2^{\mathbb{N}})^{-1}\mathbb{Z})\langle x,y\rangle$, in the other hand $G_1$ is a special Gr\"obner basis for $I$ in $R$.

\item \underline{In the ring $\mathbb{Z} _{[3]}$}, the set $G_2=\{6xyx-8xy, 4xy-3yx, 3xy^2x-3y^2x, -3yx^2+4yx, -3(yx)^2+4y^2x, \frac{64}{3}yx, y^3x, -2y^2x\}$ is a Gr\"obner basis for $(3^{\mathbb{N}})^{-1}I$ in $((3^{\mathbb{N}})^{-1}\mathbb{Z})\langle x,y\rangle$, in the other hand, $G_2$ is a special Gr\"obner basis for $I$ in $R$.
\end{itemize}
Thus the set $G=\{(G_1,(2^{\mathbb{N}})), (G_2, (3^{\mathbb{N}}))\}$ is a dynamical Gr\"obner basis for $I$ in $R$.
\end{Example}
\section{\textbf{Acknowledgment}}
Author is grateful to the IMU Berlin Einstein Foundation, the Berlin Mathematical School and to Pr.Dr. Klaus Altmann for receiving him in Berlin during the writing of the important part of this paper.


\end{document}